\documentclass[a4paper,10pt]{article}
\usepackage{amsmath}
\usepackage{amsthm}
\usepackage{amssymb}
\usepackage{latexsym}
\usepackage{enumerate}

\newtheorem{thm}{Theorem}[section]
\newtheorem{prop}{Proposition}[section]
\newtheorem{lem}{Lemma}[section]

\newtheorem{rem}{Remark}[section]

\def\bint{\displaystyle -\!\!\!\!\!\!\int}

\title{On the essential spectrum of the 
Laplacian and vague convergence of the curvature at infinity}

\author{Hironori Kumura}
\date{Department of Mathematics, Shizuoka University\\
Ohya, Shizuoka 422-8529, Japan\\
\textit{E-mail address}: smhkumu@ipc.shizuoka.ac.jp}

\begin{document}

\maketitle

\begin{abstract}
We shall prove that under some volume growth condition, 
the essential spectrum of the Laplacian contains the interval 
$[(n-1)^2K/4, \infty)$ 
if an $n$-dimensional Riemannian manifold has an end and the average 
of the part of the Ricci curvature on the end which lies below 
a nonpositive constant $(n-1)K$ 
converges to zero at infinity. 
\end{abstract}

%%%%%%    INTRODUCTION    %%%%%%

\section{Introduction}

The Laplace-Beltrami operator $\Delta$ 
on a noncompact complete Riemannian manifold $M$ is 
essentially self-adjoint on $C^{\infty}_0(M)$ 
and its self-adjoint extension to $L^2(M)$ has 
been studied by several authors from various points of view. 
For instance, Donnelly proved the following:
% 
%%%%%%    '藝 1.1    %%%%%%
%
\begin{thm}[\cite{d}]
Let $M$ be an $n$-dimensional Hadamard manifold, 
$x_0$ a fixed point of $M$, and $K$ a nonpositive constant. 
For a two plane $\pi$ in $T_xM$, 
we denote $K(x,\pi )$ to be   
the sectional curvature of $\pi$. 
Suppose that 
$$
\lim_{s\to \infty} \sup\{ |K(x,\pi )-K|;\mathrm{dist}(x,x_0)\ge s\}=0.
$$
Then the essential spectrum 
$\sigma _{\mathrm{ess}}(-\Delta)$ of $-\Delta$ 
is equal to $[-(n-1)^2K/4, \infty)$. 
\end{thm}

In this paper, we shall consider the case that the part of the 
Ricci curvature on an end of $M$ which lies below a nonpositive 
constant converges to zero in some average sense at infinity and 
some volume condition holds. And we 
study the essential spectrum of the Laplacian. 
A modification of Petersen-Wei's Lemma $2.2$ in 
\cite{p-w} is of crucial importance in our argument 
(see Lemma $2.3$ in Section $2$). 

Let us state our results more precisely. 
Let $M$ be an $n$-dimensional complete Riemannian manifold. 
Suppose that there exists an open subset $U$ of $M$ with compact smooth 
boundary $\partial U$ such that the outward pointing normal exponential 
map $\exp _{\partial U}:N^{+}(\partial U)\to M-U$ induces a diffeomorphism. 
Denote $r(*)=\mathrm{dist}(\partial U,*)$ on $M-U$ and let 
$(r,\theta )\in [0,\infty )\times \partial U$ be 
the geodesic coordinates induced from this diffeomorphism. 
We shall write the volume form as $dv_M=\sqrt{g}(r,\theta )drd\theta $ 
on the end $M-U$, where 
$d\theta$ stands for the Riemannian measure on $\partial U$ 
induced from $dv_M$. 
For latter use, 
we shall say that a plane $\pi \subset T_xM$ $(x\in M-U)$ is 
a radial plane if $\pi$ contains $\nabla r$, and by radial curvature 
we mean the restriction of the sectional curvature to all the radial 
planes. 

We introduce a new measure on the end $M-U$, which will be important tool 
in our arguments. 
For a given $C^2$-function $f:[0,\infty)\to \mathbf{R}$ satisfying 
$$\lim_{t\to \infty}f'(t)= \lim_{t\to \infty}f''(t)=0,$$ 
we introduce a measure 
$\mu_f$ on $M-U$ by 
$$\mu _f(A)=\int_A\exp (-cr(x)+f(r(x)))dv_M(x),$$
where $A$ is a measurable subset of $M-U$, $c=(n-1)\sqrt{-K}$ is a 
constant, and $K$ is a nonpositive constant. 
We also write $d\mu_f=\exp (-(n-1)\sqrt{-K}r+f(r))dv_M$. 
Moreover, we shall use the convenient notation 
$$
\frac{1}{\mu_f(A)}\int_Ah(x)d\mu_f(x)=\bint_Ah(x)d\mu_f(x),
$$
for the mean value of $h$ in $A$, when $\mu_f(A)>0$. 
This measure $d\mu_f$ is natural to consider the essential spectrum 
of the Laplacian, because  
the function $\exp (-(n-1)\sqrt{-K}r+f(r(x))$ cancels out the growth part 
of test functions and reduces our situation to subexponential growth case. 

Typical examples of the function 
$f$ are $f(t)=t^{\theta }$, $t/\log t$, $t/\log (\log t)$, and so on, 
where $\theta \in (0,1)$ is a constant. 
In this case, we should note that the weight function decomposes 
into two parts: the exponential decay part 
$\exp(-cr(x))$ and subexponential growth part $\exp(f(r))$.
The subexponential perturbation $\exp(f(r))$ will
have no effect on the essential spectrum of the Laplacian. 
 
In the sequel we will use the following notation:
$$
      A(a,b)=\{y\in M-U~|~a<\mathrm{dist}(\partial U,y)<b\}
      ~~\mathrm{for}~ 0\le a<b. 
$$
Moreover we write $h_+=\max\{h,0\}$ and 
$h_-=\max\{-h,0\}$ for a function $h$. 
In section $2$, we shall prove the following:
%
%%%%%%%    Theorem 1.2     %%%%%%%
%
\begin{thm}
Let $M$ be an $n$-dimensional complete Riemannian manifold. 
Suppose that there exists an open subset $U$ of $M$ with compact smooth 
boundary $\partial U$ such that the outward pointing normal exponential 
map $\exp _{\partial U}:N^{+}(\partial U)\to M-U$ induces a diffeomorphism. 
Denote $r(*)=\mathrm{dist}(\partial U,*)$ on $M-U$. 
We assume that the following three conditions hold:
\begin{enumerate}
	\item The measure of the end $M-U$ is infinite: 
$\displaystyle \lim_{t\to \infty}\mu_f(A(0,t))=+\infty ;$
\item there exists a constant $p>n/2$ such that \\
\hspace{20mm}
$
     \displaystyle \lim_{t\to \infty}
     \bint_{A(0,t)}\{(n-1)K-
     \mathrm{Ric}(\partial _r,\partial _r)\}_{+}^{p}d\mu_f=0;
$

\item there exists a nonnegative constant $c_1$ 
such that $\Delta r\ge -c_1$ on $M-U$. 
\end{enumerate}
Then the essential spectrum $\sigma_{\mathrm{ess}}(-\Delta )$ of $-\Delta$ 
contains the interval $[(n-1)^2K,\infty)$.
\end{thm}

We explain our assumptions in Theorem $1.2$. 
The assumption $(i)$ implies that the volume growth is greater than or 
equal to exponential except for the subexponential growth part. 
When the curvature of $M-U$ converges to a constant $K\le 0$ from {\it above} 
at infinity, $\mu_0(M-U)<\infty ~(f\equiv 0)$ may occurs. 
But choosing $f$ appropriately, 
we will be able to get the situation $\mu_f(M-U)=\infty$ in many cases. 
The assumption $(i)$ and $(ii)$ mean that the average of the part of the 
Ricci curvature on the end $M-U$ which lies below 
a nonpositive constant $(n-1)K$ converges to zero {\it at infinity}. 

Next we should add some comments about the condition $(iii)~\Delta r\ge -c_1$. 
If $M$ satisfies the assumptions of Theorem $1.1$, 
$\Delta r$ approaches to $(n-1)\sqrt{-K}$ at infinity 
and hence  
the condition $(iii)$ $\Delta r \ge -c_1$ in Theorem $1.2$ is 
indeed satisfied. 
We also note the following:
%%
%%%%%   Proposition   %%%%%
%%
\begin{prop}
Let $M$ be as in Theorem $1.2$. 
\begin{enumerate}
\item If the mean curvature of $\partial U$ is nonpositive, i.e., 
$\Delta r\ge 0$ on $\partial U$ and radial curvature on $M-U$ is 
nonpositive, then $\Delta r\ge 0$ on $M-U$.
\item If there exists a positive constant $A(>K)$ such that 
$\mathrm{Ric}(\partial r,\partial r)\ge -(n-1)A$ on $M-U$, 
then $\Delta r\ge -(n-1)\sqrt{A}$ on $M-U$.
\end{enumerate}
\end{prop} 
Thus, if either of the conditions in Proposition $1.1$ holds, then 
the condition $(iii)$ $\Delta r \ge -c_1$ in Theorem $1.2$ is satisfied. 
Proposition $1.1$ is proved by using the 
comparison theorem in Riemannian geometry in section $2$. 

The final remark is that 
in Theorem $1.2$ we do not assume that the average of 
the part of the curvature 
which lies {\it above} $K$ converges to zero at infinity, 
but we only assume that 
$(i)~\displaystyle \lim_{t\to \infty}\mu_f(A(0,t))=+\infty $. 
This condition $(i)$ bounds the curvature from the above very mildly,  
which seems to be worthy of note. 
%
%%%%%%    SECTION 2    %%%%%%
%
\section{Proof}
In this section, we shall prove Theorem $1.2$ by transplantation method 
and Proposition $1.1$ by the comparison theorem in Riemannian geometry. 

Theorem $1.2$ readily follows from the following five lemmas. 
In the following, we simply denote $c=(n-1)\sqrt{-K}$. 
First, we note the following obvious lemma:
%%%%%%    LEMMA 2.1    %%%%%%
%
\begin{lem}
Let $M$ be an $n$-dimensional complete Riemannian manifold. 
Suppose that there exists an open subset $U$ of $M$ with compact smooth 
boundary $\partial U$ such that the outward pointing normal exponential 
map $\exp _{\partial U}:N^{+}(\partial U)\to M-U$ induces a diffeomorphism. 
For any constants $\delta \in (0,1)$ 
and monotone increasing sequence $\{b_i\}$ 
with $\displaystyle \lim_{i\to \infty}b_i=\infty$, 
define the sequence $\{a_i\}$ by 
$\mu_f(A(0,a_i))=\delta\cdot \mu_f(A(0,b_i))$. 
Then the following two conditions are equivalent:\vspace{1mm}\\
\hspace{2mm}$(i)$~$\displaystyle\lim_{i\to \infty}a_i=\infty$.\vspace{1mm}\\
\hspace{2mm}$(i)$~$\displaystyle\lim_{t\to \infty}\mu_f(A(0,t))=\infty.$
\end{lem}
%
%%%%%%    LEMMA 2.2    %%%%%%
%
\begin{lem}
Let $M$ be as above. 
Suppose that for any constant $\delta \in (0,1)$ 
and monotone increasing sequence $\{d_i\}$ 
with $\displaystyle\lim_{i\to \infty}d_i=\infty$, 
there exists a sequence $\{a_i\}$ such that\vspace{1mm}\\
\hspace{2mm}$(i)$~$\mu_f(A(0,a_i))=\delta\cdot \mu_f(A(0,d_i))$;\vspace{1mm}\\
\hspace{2mm}$(ii)$~$\displaystyle\lim_{i\to \infty}a_i
=\lim_{i\to \infty}(d_i-a_i)=\infty$;\vspace{1mm}\\
\hspace{2mm}$(iii)$~$\displaystyle\lim_{i\to \infty}\bint_{A(a_i,d_i)}
|\Delta r-c|^2d\mu_f=0$.\vspace{1mm}\\
Then the essential spectrum of $-\Delta$ contains the 
interval $[c^2/4,\infty)$. 
\end{lem}
\begin{proof}
For an arbitrary increasing sequence $\{d_i\}$ 
with $\displaystyle \lim_{i\to \infty}d_i=\infty$, we define 
sequences $0<a_i<b_i<c_i<d_i$ $(i=1,2,\dots )$ 
by the following equation:
$$
   \mu_f(A(0,a_i))=\mu_f(A(a_i,b_i))
   =\mu_f(A(b_i,c_i))=\mu_f(A(c_i,d_i)).
$$
Then, our assumptions imply
\begin{align}
    &\lim_{i\to \infty}a_i
   =\lim_{i\to \infty}(b_i-a_i)
   =\lim_{i\to \infty}(c_i-b_i)
   =\lim_{i\to \infty}(d_i-c_i)
  =\infty ,\\
  &\lim_{i\to \infty}\bint_{A(a_i,d_i)}
  |\Delta r-c|^2d\mu_f=0.
\end{align}
Indeed, first, we see that $(i)$ and $(ii)$ imply that 
$\displaystyle\lim_{i\to \infty}a_i=\lim_{i\to \infty}b_i
=\lim_{i\to \infty}c_i=\lim_{i\to \infty}(d_i-c_i)=\infty$. 
Then, since $\frac{1}{2}\mu_f(A(0,b_i))=\mu_f(A(0,a_i))$ and 
$\displaystyle\lim_{i\to \infty}b_i=\infty$, 
$(i)$ and $(ii)$ with $d_i$ replaced by $b_i$ 
imply $\displaystyle\lim_{i\to \infty}(b_i-a_i)=\infty$. 
Similarly, we see that $\displaystyle\lim_{i\to \infty}(c_i-b_i)=\infty$. 

In view of (1), we can choose a sequence of functions 
$\varphi_i\in C^2_0(\mathbf{R})$ such that 
\begin{align*}
   &\varphi_i(t)=1 \quad (b_i\le t \le c_i); 
   ~~\varphi_i(t)=0 \quad (t\le a_i,~~d_i\le t);\\ 
   &0\le \varphi_i(t)\le 1~~\mathrm{for ~all}~ t\in \mathbf{R};\\
   &\sup |\varphi '|\le A_i,
   ~~\sup |\varphi ''|\le A_i,~~\lim_{i\to \infty}A_i=0.
\end{align*}
For any $\lambda >c^2/4$, we define functions 
$$
   \phi _i(x):=\sqrt{v_i}\cdot \varphi_i(r(x))\cdot 
   \exp\left(ar(x)+\frac{f(r(x))}{2}\right)
    \in C^2_0(M-U),
$$
where we set $v_i=\mu_f(A(a_i,d_i))^{-1}$ and 
$\displaystyle a=-\frac{c}{2}+\sqrt{-\left(\lambda-\frac{c^2}{4}\right)}
\in \mathbf{C}$. 
Then, a direct computation shows that 
\begin{align*}
   (\Delta \phi_i+\lambda \phi_i)
   =&\sqrt{v_i}\varepsilon(r)
     \exp \left(ar+\frac{f(r)}{2}\right)\\
   &  +\sqrt{v_i}\{\varphi_i'(r)+\varphi_i(r)f'(r)/2\}
     \cdot \Delta r\cdot \exp \left(ar+\frac{f(r)}{2}\right)\\
   &+\sqrt{v_i} a\varphi_i(r) (\Delta r-c) 
      \cdot \exp \left(ar+\frac{f(r)}{2}\right),
\end{align*}
where we set 
\begin{align*}
     \varepsilon(r)=&\varphi_i''(r)+2a \varphi_i'(r)+\varphi_i'(r)f'(r)
                      +\varphi_i(r)f''(r)/2\\
      &+a\varphi_i(r)f'(r)+\varphi_i(r)(f'(r))^2.
\end{align*}
Hence, we have 
\begin{align*}
   \int_{M-U}
   |\Delta \phi_i+\lambda \phi_i|^2dv_M
   \le &c_2\left(A_i^2+A_i^4+B_i^2+B_i^4\right)
   \,\bint_{A(a_i,d_i)}d\mu_f\\
   &+2\left(A_i^2+B_i^2\right)
   \,\bint_{A(a_i,d_i)}|\Delta r|^2 d\mu_f\\
   &+2|a|^2
   \,\bint_{A(a_i,d_i)}|\Delta r-c|^2d\mu_f,
\end{align*}
where $\displaystyle B_i=\max\{\sup_{t\ge a_i} |f'(t)|,\sup_{t\ge a_i} |f''(t)|\}$ 
and $c_2$ is a constant which depends only on $a$. 
Since 
$$
    A_i^2\bint_{A(a_i,d_i)}|\Delta r|^2d\mu_c 
    \le 2A_i^2\left(|c|^2+
    \bint_{A(a_i,d_i)}|\Delta r-c|^2d\mu_f\right),
$$
we see from $(2)$ and $\displaystyle\lim_{i\to \infty}B_i=0$ that 
\begin{equation}
     \lim_{i\to \infty }\int_M|\Delta \phi_i+\lambda \phi_i|^2dv_M=0.
\end{equation}
On the other hand, 
\begin{equation}
   1\ge \int_M|\phi_i|^2dv_M
   \ge v_i\cdot \mu_f(A(b_i,c_i))=\frac{1}{3}.
\end{equation}
In view of our construction of $\{\phi_i\}$, $(3)$ and $(4)$ imply that 
$\lambda \in \sigma_{\mathrm{ess}}(-\Delta)$. 
Since 
$\lambda > c^2/4$ is arbitrary, 
we obtain $[c^2/4,\infty)\subset \sigma_{\mathrm{ess}}(-\Delta)$. 
\end{proof}

A modification of the argument of the proof of 
Lemma $2.2$ in Petersen-Wei \cite{p-w} 
yields the following lemma, which is of crucial importance in our argument.
%
%%%%%%    LEMMA 2.3    %%%%%%
%
\begin{lem}
Let $M$ be an $n$-dimensional complete Riemannian manifold. 
Suppose that there exists an open subset $U$ of $M$ with compact smooth 
boundary $\partial U$ such that the outward pointing normal exponential 
map $\exp _{\partial U}:N^{+}(\partial U)\to M-U$ induces a diffeomorphism. 
Let $K \le 0$ be a constant 
and set $c=(n-1)\sqrt{-K}$. 
As a technical condition, we assume that 
$c\ge f'(t)$ for $t\ge 0$. 
Then if $p>n/2$, we have 
\begin{align*}
   &\int^r_0\{\Delta r(r,\theta)
   -h_K(r)\}_{+}^{2p}
      \exp(-cr+f(r))\sqrt{g}(r,\theta)dr\nonumber \\
   &\le c_2(n,p)\int^r_0\{(n-1)K
      -\mathrm{Ric}(\partial_r,\partial_r)\}_{+}^p
      \exp(-cr+f(r))\sqrt{g}(r,\theta)dr,
\end{align*}
where $c_2(n,p)$ is a constant depending only on $n$ and $p$, 
and $h_K$ is the solution to 
\begin{equation}
    h_K'(t)+\frac{h_K^2(t)}{n-1}=-(n-1)K,
\end{equation}
with initial condition $h_K(0)=a_0$, 
where $a_0\ge 0$ is a constant satisfying 
$\displaystyle a_0\ge \max_{\partial U} \Delta r$. 
\end{lem}
\begin{proof}
We set $\Delta r=h(r,\theta )$, 
then $h$ satisfies the differential inequality
$$
     \partial _r h+\frac{h^2}{n-1}\le -\mathrm{Ric}(\partial_r,\partial_r).
$$
For simplicity, we define
\begin{align*}
&\rho(r,\theta):=\max\{0,(n-1)K-\mathrm{Ric}(\partial_r,\partial_r)\},\\
&\psi (r,\theta):=\max\{0,h(r,\theta)-h_K(r)\}
\end{align*}
and set $\widetilde{\sqrt{g}}(r,\theta):=\exp(-cr+f(r))\sqrt{g}(r,\theta)$.
In the following, we fix $\theta \in \partial U$, therefore we simply 
write $\rho(r)$ and $\psi (r)$ for $\rho(r,\theta)$ and $\psi (r,\theta)$ 
respectively. 
Then we see that $\psi$ is absolutely continuous and satisfies 
\begin{equation}
     \psi '+\frac{\psi^2}{n-1}+2\frac{\psi h_K}{n-1}\le \rho.
\end{equation}
Indeed, if $(h-h_K)(r_0)>0$, then $\psi(r)$ coincides with 
$(h-h_K)(r)$ around $r_0$. 
Hence
\begin{align*}
   &\psi'+\frac{\psi^2}{n-1}+2\frac{\psi h_K}{n-1}
   =h'-h'_K+\frac{h^2}{n-1}-\frac{h_K^{2}}{n-1}\\
   \le &(n-1)K-\mathrm{Ric}(\partial r,\partial r)\le \rho
\end{align*}
and $(6)$ holds around $r_0$. 
When $(h-h_K)(r_0)=0$, $\psi(r_0)=0$ by definition. 
Moreover, in this case, we may take $\psi'(r_0)=0$ and 
$(6)$ also holds at these points because $\rho$ is nonnegative. 
When $(h-h_K)(r_0)<0$, $\psi(r)=0$ around $r_0$ and hence 
the left hand side of $(6)$ is zero. 
Therefore $(6)$ also holds at this point. 
Thus we conclude that $(6)$ holds for all $r$. 

Multiplying $(6)$ by $\psi^{2p-2}\exp\{-cr+f(r)\}\sqrt{g}(r,\theta)$ 
and integrating over $[0,r]$, we get 
\begin{align*}
   \int_0^r\rho\psi^{2p-2}\widetilde{\sqrt{g}}(t,\theta)dt
   \ge &\int_0^r\psi'\psi^{2p-2}\widetilde{\sqrt{g}}(t,\theta)dt
    +\frac{1}{n-1}\int_0^r\psi^{2p}\widetilde{\sqrt{g}}(t,\theta)dt\\
     &+\frac{2}{n-1}\int_0^rh_K\psi^{2p-1}\widetilde{\sqrt{g}}(t,\theta)dt.\\
\end{align*}
Integration by parts yields
\begin{align*}
   &\int_0^r\psi'\psi^{2p-2}\widetilde{\sqrt{g}}(t,\theta)dt\\
   =&\left[\frac{\psi^{2p-1}\widetilde{\sqrt{g}}(t,\theta)}
   {2p-1}\right]^r_0 -\frac{1}{2p-1}\int_0^r\psi^{2p-1}\{-c+f'(r)+h\}
   \widetilde{\sqrt{g}}(t,\theta)dt\\
   \ge &-\frac{1}{2p-1}\int^r_0\psi^{2p-1}(t)(f'(t)-c+h_K(t))
   \widetilde{\sqrt{g}}(t,\theta)dt\\
   &-\frac{1}{2p-1}\int_0^r\psi^{2p}\widetilde{\sqrt{g}}(t,\theta)dt,
\end{align*}
where we have used $\psi(0)=0$, which follows from our assumption 
$h_K(0)=a_0\ge \max \{\Delta r(0,\theta)~|~\theta\in \partial U\}$. 
Inserting this inequality in the previous one, we have
\begin{align*}
     \int_0^r\rho(t)\psi^{2p-2}(t)\widetilde{\sqrt{g}}(t,\theta)dt
\ge & \int_0^r\psi^{2p-1}(t)\frac{w(t)}{(n-1)(2p-1)}
       \widetilde{\sqrt{g}}(t,\theta)dt\\
    &  +\left(\frac{1}{n-1}-\frac{1}{2p-1}\right)\int_0^r\psi^{2p}(t)
       \widetilde{\sqrt{g}}(t,\theta)dt.
\end{align*}
Here, we set 
$w(t)=(n-1)(c-f'(t))+(4p-n-1)h_K(t)$. 
Since $h_K(t)$ satisfies $(5)$ and the initial condition $h_K(0)=a_0\ge 0$, 
we see that $h_K(t)\ge 0$ for $t\ge 0$ and, 
hence, $w(t)\ge 0$ for all $t\ge 0$. 
Thus we obtain 
\begin{align*}
&\left(\frac{1}{n-1}-\frac{1}{2p-1}\right)
\int_0^r\psi^{2p}(t)\widetilde{\sqrt{g}}(t,\theta)dt\\
&\le \int_0^r\rho(t)\psi^{2p-2}(t)\widetilde{\sqrt{g}}(t,\theta)dt\\
&\le \left(\int^r_0\rho^p(t)\widetilde{\sqrt{g}}(t,\theta)dt\right)^{1/p}
\cdot \left(\int^r_0\psi^{2p}(t)\widetilde{\sqrt{g}}(t,\theta)dt\right)^{1-\frac{1}{p}}.
\end{align*}
Hence we have 
$$
\left(\frac{1}{n-1}-\frac{1}{2p-1}\right)
\left(\int^r_0\psi^{2p}(t)\widetilde{\sqrt{g}}(t,\theta)dt\right)^{1/p}
\le \left(\int^r_0\rho^p(t)\widetilde{\sqrt{g}}(t,\theta)dt\right)^{1/p},
$$
that is, 
$$
\int^r_0\psi^{2p}(t)\widetilde{\sqrt{g}}(t,\theta)dt
\le \left(\frac{1}{n-1}-\frac{1}{2p-1}\right)^{-p}
\int^r_0\rho^p(t)\widetilde{\sqrt{g}}(t,\theta)dt.
$$
Thus, we obtain the desired inequality with 
$c_2(n,p)=\left(\frac{1}{n-1}-\frac{1}{2p-1}\right)^{-p}$.
\end {proof}
%
%%%%%%    REMARK 2.1    %%%%%%
%
\begin{rem}
We need not worry about 
the technical condition $c\ge f'(t)$ in Lemma $2.3$,
because our condition 
$\displaystyle \lim_{t\to \infty}\mu_f(A(0,t))=+\infty$ 
will enables us to disperse it. 
\end{rem}
%
%%%%%%    LEMMA 2.4    %%%%%%
%
\begin{lem}
Let $M,r,K,c=(n-1)\sqrt{-K}$ be as above and assume the 
following:\vspace{1mm}\\
\hspace{2mm}$(i)$~
$\displaystyle\lim_{t\to \infty}\mu_f(A(0,t))=+\infty ;$\vspace{1mm}\\
\hspace{2mm}$(ii)$~There exists a constant $p>n/2$ such that \vspace{1mm}\\
\hspace{20mm}
$
     \displaystyle \lim_{t\to \infty}
     \bint_{A(0,t)}\{(n-1)K-
     \mathrm{Ric}(\partial _r,\partial _r)\}_{+}^{p}d\mu_f=0.
$\\
Then we have
\begin{equation}
     \lim_{t\to \infty}\frac{d}{dt}\log (\mu_f(A(0,t))=0.
\end{equation}
In particular, for any monotone increasing sequence $\{b_i\}$ with 
$\displaystyle \lim_{i\to \infty}b_i=\infty$ and $\delta\in (0,1)$, 
define a sequence 
$\{a_i\}$ by $\mu_f(A(a_i,b_i))=\delta \cdot \mu_f(A(0,b_i))$. 
Then we have $\displaystyle \lim_{i\to \infty}(b_i-a_i)=\infty$. 
\end{lem}
\begin{proof}
We again set 
$\widetilde{\sqrt{g}}(r,\theta):=\exp(-cr+f(r))\sqrt{g}(r,\theta)$ 
for the sake of simplicity. 
First we note that $\displaystyle \lim_{r\to \infty}h_K(r)=(n-1)\sqrt{-K}$ 
because $h_K(0)=a_0\ge 0$.  
Let $\varepsilon >0$ be given and take sufficiently large $s_1>0$. 
For $t>s_1$, integrating the equation 
$\partial _r\widetilde{\sqrt{g}}(r,\theta)
=\{(\Delta r)(r,\theta)-c+f'(r)\}\widetilde{\sqrt{g}}(r,\theta)$ 
over $[s_1,t]\times \partial U$, we have 
\begin{align*}
     &  \int_{\partial U}\widetilde{\sqrt{g}}(t,\theta)d\theta
        -\int_{\partial U}\widetilde{\sqrt{g}}(s_1,\theta)d\theta\\
    =&  \int^t_{s_1}dr\int_{\partial U}
        \{\Delta r-c+f'(r)\}
        \widetilde{\sqrt{g}}(r,\theta)d\theta\\
    =&  \int^t_{s_1}dr\int_{\partial U}
        \{\Delta r-h_K(r)\}
        \widetilde{\sqrt{g}}(r,\theta)d\theta\\
     &  +\int^t_{s_1}dr\int_{\partial U}
        \left(h_K(r)-c+f'(r)\right)
        \widetilde{\sqrt{g}}(r,\theta)d\theta\\
  \le&  \int_{A(0,t)}
        \left\{\Delta r-h_K(r)\right\}_{+}d\mu_f
        +\varepsilon \mu_f(A(s_1,t))\\
  \le&  \mu_f(A(0,t))\left(\bint_{A(0,t)}
        \left\{\Delta r-h_K(r)\right\}_{+}^{2p}
        d\mu_f\right)^{\frac{1}{2p}}+\varepsilon \mu_f(A(s_1,t))\\
  \le&  \mu_f(A(0,t))c_3(n,p)\left(\bint_{A(0,t)}
        \left\{(n-1)K-\mathrm{Ric}(\partial_r,\partial_r)\right\}_{+}^{p}
        d\mu_f\right)^{\frac{1}{2p}}\\
     &  +\varepsilon \mu_f(A(s_1,t)),
\end{align*}
In the last inequality above, we have used Lemma $2.3$. 
Hence, 
\begin{align*}
     (0\le &)\frac{d}{dt}\log \mu_f\left(A(0,t)\right)
             =\mu_f(A(0,t))^{-1}
             \int_{S^{n-1}(1)}\widetilde{\sqrt{g}}(t,\theta)d\theta\\
     \le &\mu_f(A(0,t))^{-1}
     \int_{S^{n-1}(1)}\widetilde{\sqrt{g}}(s_1,\theta)d\theta\\
     &+c_3(n,p)\left(\bint_{A(0,t)}
       \{(n-1)K-\mathrm{Ric}(\partial_r,\partial_r)\}_{+}^{p}
       d\mu_f\right)^{\frac{1}{2p}}+\varepsilon
\end{align*}
and our assumptions $(i)$ and $(ii)$ imply that the equation $(7)$ holds. 

Next we shall show that $(7)$ implies 
$\displaystyle\lim_{i\to \infty}(b_i-a_i)=\infty$. 
For any $\varepsilon >0$, there exists $t_0>0$ such that 
\begin{equation}
     \frac{d}{dt}\log (\mu_f(A(0,t)))\le \varepsilon \qquad 
     \mathrm{for~all}~t\ge t_0.
\end{equation}
Since $\displaystyle\lim_{i\to \infty}a_i=\infty$ by Lemma $2.1$, 
there exits $i_0$ such that $a_i>t_0$ for all $i\ge i_0$. 
Thus, for any $i\ge i_0$, integrating this inequality $(7)$ over $[a_i,b_i]$, 
we obtain
$$
     0<-\log(1-\delta)
     =\log\frac{\mu_f(A(0,b_i))}{\mu_f(A(0,a_i))}
     \le \varepsilon(b_i-a_i).
$$
Hence, $\displaystyle\lim_{i\to \infty}(b_i-a_i)=\infty$. 
\end{proof}
%
%%%%%%    LEMMA 5    %%%%%%
%
\begin{lem}
We assume that the following three conditions hold:\vspace{1mm}\\
\hspace{2mm}$(i)$~$\displaystyle\lim_{t\to \infty}\mu_f(A(0,t))=+\infty ;$
\vspace{1mm}\\
\hspace{2mm}$(ii)$~there exists a constant $p>n/2$ such that \vspace{1mm}\\
\hspace{20mm}
$
     \displaystyle \lim_{t\to \infty}
     \bint_{A(0,t)}\{(n-1)K-
     \mathrm{Ric}(\partial _r,\partial _r)\}_{+}^{p}d\mu_f=0;
$\vspace{1mm}\\
\hspace{2mm}$(iii)$~there exists a constant 
$c_1\ge 0$ such that $\Delta r\ge -c_1$. \vspace{1mm}\\
If we define the sequence $\{a_i\}$ by 
$\mu_f(A(0,a_i))=\delta\cdot \mu_f(A(0,b_i))$ for any monotone 
increasing sequence $\{b_i\}$ 
with $\displaystyle\lim_{i\to \infty}b_i=\infty$ and $\delta \in (0,1)$, 
then the condition $(iii)$ in Lemma $2.2$ holds. 
\end{lem}
\begin{proof}
By Lemma $2.3$, we have 
\begin{align*}
       & \mu_f(A(a_i,b_i))^{-1}\int_{\partial U}d\theta 
         \int^{b_i}_{a_i}\{(\Delta r)(r,\theta)-h_K(r)\}^{2p}_{+}
         \widetilde{\sqrt{g}}(r,\theta)dr\\
   \le & \frac{c_2(n,p)}{(1-\delta)\mu_f(A(0,b_i))}
          \int_{\partial U}d\theta \int^{b_i}_{0}
          \{(n-1)K-\mathrm{Ric}(\partial _r,\partial _r)\}^p_{+}
          \widetilde{\sqrt{g}}(r,\theta)dr\\
     &   \to 0\qquad \mathrm{as}\quad i\to \infty.
\end{align*}
Hence, in view of the fact that 
$\displaystyle\lim_{i\to \infty}a_i=\infty$ and 
$\displaystyle\lim_{r\to \infty}h_K(r)=(n-1)\sqrt{-K}=c$, we see that
\begin{align}
     &  \bint_{A(a_i,b_i)}
       \{\Delta r-c\}_{+}^2d\mu_f\nonumber \\ 
  \le& \left(\bint_{A(a_i,b_i)}\{\Delta r-c\}_{+}^{2p}d\mu_f\right)
       ^{\frac{1}{p}}
       \to 0\qquad \mathrm{as}\quad i\to \infty.
\end{align}
On the other hand, integrating the equation 
$\displaystyle
     \widetilde{\sqrt{g}}(b_i,\theta)-\widetilde{\sqrt{g}}(a_i,\theta)
     =\int_{a_i}^{b_i}\partial _r\widetilde{\sqrt{g}}(r,\theta)dr
$ over $\partial U$, we have 
\begin{align*}
     &\mu_f(A(a_i,b_i))^{-1}\int_{\partial U}\widetilde{\sqrt{g}}
     (b_i,\theta)d\theta
     -\mu_f(A(a_i,b_i))^{-1}\int_{\partial U}
     \widetilde{\sqrt{g}}(a_i,\theta)d\theta\\
     =&\mu_f(A(a_i,b_i))^{-1}\int_{\partial U}d\theta
     \int_{a_i}^{b_i}\partial _r\widetilde{\sqrt{g}}(r,\theta)dr.
\end{align*}
Here, since
\begin{align*}
     &\mu_f(A(a_i,b_i))^{-1}\int_{\partial U}\widetilde{\sqrt{g}}
     (b_i,\theta)d\theta\\
     =&\frac{\mu_f(A(0,b_i))}{\mu_f(A(a_i,b_i))}
     \cdot \mu_f(A(0,b_i))^{-1}\int_{\partial U}\widetilde{\sqrt{g}}
     (b_i,\theta)d\theta\\
     =&(1-\delta)^{-1}\frac{d}{dt}\Big|_{t=b_i}\log\mu_f(A(0,t))
     \to 0\qquad \mathrm{as}\quad i\to \infty
\end{align*}
and similarly 
\begin{align*}
     &  \mu_f(A(a_i,b_i))^{-1}\int_{\partial U}
       \widetilde{\sqrt{g}}(a_i,\theta)d\theta\\
    =&\frac{\delta}{1-\delta}
     \cdot \frac{d}{dt}\Big|_{t=a_i}\log\mu_f(A(0,t))
     \to 0\qquad \mathrm{as}\quad i\to \infty,
\end{align*}
we have 
$$
     \lim_{i\to \infty}
       \bint_{A(a_i,b_i)}\{\Delta r-c\}d\mu_f
    =\lim_{i\to \infty}
       \bint_{A(a_i,b_i)}\{\Delta r-c+f'(r)\}d\mu_f=0.
$$
Therefore, in view of H\"older's inequality, $(9)$ implies that 
$$
     \lim_{i\to \infty}
     \bint_{A(a_i,b_i)}\{\Delta r-c\}_{-}d\mu_f=0.
$$
Here, since
$
     \{\Delta r-c\}_{-}=\max\{0,-(\Delta r-c)\}\le c_1+c
$
on $M-U$, 
\begin{equation}
     \lim_{i\to \infty}
     \bint_{A(a_i,b_i)}\{\Delta r-c\}_{-}^2d\mu_f=0.
\end{equation}
The equations $(9)$ and $(10)$ give the desired equality. 
\end{proof}
From five lemmas above and Remark $2.1$, we see that Theorem $1.2$ holds. 

\vspace{5mm}
\noindent{\it Proof of Proposition $1.1$}\vspace{1mm}\\
Since $(i)$ immediately follows from the comparison theorem in Riemannian 
geometry, we shall omit its proof.
Hence we shall only prove $(ii)$. 
Assuming that there exists a point $(r_0,\theta_0)\in M-U$ such that 
$(\Delta r)(r_0,\theta_0)<-(n-1)\sqrt{A}$, we shall obtain a contradiction. 
Let us set $h(r)=(\Delta r)(r,\theta_0)$ for $r\ge r_0$. 
Moreover let $f(r):[r_0,l)\to \mathrm{R}$ be the 
solution to the differential equation
$$
   f'+\frac{f^2}{n-1}-(n-1)A=0;~~f(r_0)=(\Delta r)(r_0,\theta_0),
$$
where $[r_0,l)$ is the maximal interval of existence of the solution. 
Then the comparison theorem in Riemannian geometry implies 
$h(r)\le f(r)$ on $[r_0,l)$. 
But the solution $f$ is expressed as follows:
$$
   f(r)=-(n-1)\sqrt{A}\frac{1+a_0e^{2\sqrt{A}r}}{1-a_0e^{2\sqrt{A}r}},
$$
where $a_0>0$ is the constant determined by 
$a_0e^{2\sqrt{A}r_0}=(B-1)/(B+1)$ and 
$B=-(\Delta r)(r_0,\theta_0)/(n-1)\sqrt{A}>1$. 
Hence 
$\displaystyle 
l=r_0+\frac{1}{2\sqrt{A}}\log \left(\frac{B+1}{B-1}\right)<\infty$ 
and $\displaystyle\lim_{r\to l-0}f(r)=-\infty$. 
This contradicts the fact that 
$h(r)$ is defined for all $r\ge r_0$. 
Therefore, $\Delta r\ge -(n-1)\sqrt{A}$ on $M-U$ and 
this completes the proof of Proposition $1.1$. \vspace{5mm}

Finally, we note that 
the condition $(i)$ can be replaced by a more familiar one for 
geometers:
\begin{lem}
Let us set $V(t)=\mathrm{Vol}(A(0,t))$. 
Then we see the following: 
\begin{enumerate}
	\item If $c>0$, 
$\displaystyle\lim_{t\to \infty}\exp(-ct+f(t))V(t)=\infty$ implies  
$\displaystyle\lim_{t\to \infty}\mu_f(A(0,t))=+\infty$. 
	\item If $c=0$, we add an assumption $f'\ge0$. 
Then $\displaystyle\lim_{t\to \infty}\exp(f(t))V(t)=\infty$ implies  
$\displaystyle\lim_{t\to \infty}\mu_f(A(0,t))=+\infty$. 
\end{enumerate}
\end{lem}
\begin{proof}
We write $S(r)=\mathrm{Vol}_{n-1}(\partial_{+} A(0,r))$, where 
$\partial_{+} A(0,r)=\{y\in M-U | $\\
$\mathrm{dist}(\partial U, y)=r\}$. 
Then integration by parts yields
\begin{align*}
     \mu_f(A(0,t))&=\int_0^t\exp(-cr+f(r))S(r)dr\\
     &=\left[\exp(-cr+f(r))V(r)\right]^t_0
     -\int_0^t(-c+f'(r))\exp(-cr+f(r))V(r)dr\\
     &=\exp(-ct+f(t))V(t)+\int_0^t(c-f'(r))\exp(-cr+f(r))V(r)dr.
\end{align*}
In view of our assumption $\displaystyle \lim_{t\to \infty}f'(t)=0$, 
we immediately see that Lemma $2.6$ holds. 
\end{proof}

%%%%%%    ŽQl•¶Œ£    %%%%%%


\begin{thebibliography}{99}
\bibitem{b1}~R.~Brooks, 
\textit{A relation between growth and the spectrum of the Laplacian},
 Math.~Z. \textbf{178}~(1981), 501-508.
%
\bibitem{b2}~R.~Brooks, 
\textit{On the spectrum of non-compact manifolds with finite volume},
 Math.~Z., \textbf{187}~(1984), 425-432.
%
\bibitem{c}~I.~Chavel, 
{\it Eigenvalues in Riemannian Geometry}, Academic Press, New York, 1984.

%
\bibitem{d}~H.~Donnelly,~
{\it On the essential spectrum of a complete Riemannian manifold},
 Topology~\textbf{20}~(1981), 1-14.
%
\bibitem{d-g}~H.~Donnelly~and~N.~Garofalo, 
{\it Riemannian manifolds whose Laplacians have purely continuous spectrum},
 Math.~Ann., \textbf{293}~(1992), 143-161.
%
\bibitem{d-l}~H.~Donnelly~and~P.~Li, 
{\it Pure point spectrum and negative curvature for noncompact manifolds},
 Duke~Math.~J.,~\textbf{46}~(1979),~497-503.
%
\bibitem{e-w}~K.~D.~Elworthy and F.-Y.~Wang, 
{\it On the essential spectrum of the Laplacian on Riemannian manifolds},
 preprint.
%
\bibitem{e}~J.~F.~Escobar, 
{\it On the spectrum of the Laplacian on complete Riemannian manifolds},
 Comm.~Partial~Differential Equations, \textbf{11}~(1986), 63-85.
%
\bibitem{e-f}~J.~F.~Escobar~and A.~Freire,
{\it The spectrum of the Laplacian of manifolds of positive curvature},
~Duke~Math.~J., \textbf{65}~(1992), 1-21.
%
\bibitem{g}~S.~Gallot, 
{\it Isoperimetric inequalities based on integral norms of Ricci curvature},
 Ast\'{e}risque 157-158~(1988), 191-216.
%
\bibitem{k-l}~L.~Karp, 
{\it Noncompact manifolds with purely continuous spectrum},
 Michigan~Math.~J., \textbf{31}~(1984), 339-347.
%
\bibitem{k}~H.~Kumura, 
{\it On the essential spectrum of the Laplacian on complete manifolds},
 J.~Math.~Soc.~Japan, \textbf{49}~(1997), 1-14.
%
\bibitem{l}~J.~Li,~
{\it Spectrum of the Laplacian on a complete Riemannian 
manifold 
with nonnegative Ricci curvature which possess(es) a pole},
 J.~Math.~Soc.~Japan, \textbf{46}~(1994), 213-216.
%
\bibitem{p-w}~P.~Petersen and G.~Wei, 
\textit{Relative volume comparison with integral curvature bounds},
 GAFA, Geom.~funct.~anal., \textbf{7}~(1997), 1031-1045.
%
\end{thebibliography}
\end{document}